\documentclass[10pt]{amsart}
    \usepackage[all]{xy}
    \usepackage{amssymb}
    \usepackage{fancybox}
    \usepackage{oldgerm}
\usepackage{graphicx}

    \newcommand{\Spe}{{\rm Spec}}

    \newcommand{\be}{\begin{equation}}
    \newcommand{\ee}{\end{equation}}

    \newtheorem{pro}{Proposition}[section]

    \newtheorem{cor}[pro]{Corollary}
    \newtheorem{theorem}[pro]{Theorem}

\begin{document}
    \bibliographystyle{amsplain}

\title{On the principal ideal theorem in arithmetic topology}
\author{Dimoklis Goundaroulis  \and Aristides Kontogeorgis }

\address{
National Technical University of Athens, Department of Mathematics\\
} \email{dgound@gmail.com}

\address{
Department of Mathematics, University of the \AE gean, 83200 Karlovassi, Samos,based  on the automorphism groups
Greece\\ { \texttt{\upshape http://eloris.samos.aegean.gr}}
}
\email{kontogar@aegean.gr}

\date{\today}

\maketitle

    \begin{abstract}
In this paper we state and prove the analogous of the principal ideal theorem of algebraic number theory
for the case of $3$-manifolds from the point of view of arithmetic topology.
        \end{abstract}


    \maketitle

\section{Introduction}

There are
certain analogies between the notions  of number theory and those of
3-dimensional topology, that are described by the following
dictionary, named after Mazur, Kapranov and Reznikov.
\begin{itemize}

\item Closed, oriented, connected, smooth 3-manifolds correspond
to affine  schemes $\Spe \mathcal{O}_{K}$, where K is an algebraic
number field and $\mathcal{O}_K$ denotes the ring of algebraic integers of $K$.

\item A link in $M$ corresponds to an ideal in $\mathcal{O}_K$ and
a knot in $M$ corresponds to a prime ideal in $\mathcal{O}_K$.

\item An algebraic integer $w \in \mathcal{O}_K$ is analogous to
an embedded surface (possibly with boundary).


\item The class group $\mathrm{Cl}(K)$ corresponds to $H_1(M,\mathbb{Z})$.

\item Finite extensions of number fields  $L/K$ correspond to
finite branched coverings of 3-manifolds $\pi: M \rightarrow N$. A
\textit{branched cover} $M$ of a 3-manifold $N$ is given by a map
$\pi$  such that there is a link $L$ of $N$ with the following
property: The restriction map  $\pi:M \backslash \pi^{-1}(L)
\rightarrow N \backslash L $ is a topological cover.

\end{itemize}

For the necessary background in algebraic number theory the reader
should look at any standard book, for example \cite{Janusz}. For
the topological part: by the term {\em knot} (resp. {\em link}) we
mean  {\em tame knot} (resp. {\em tame link}). By the term  {\em
embedded surface} we mean   an embedding  $f:E\rightarrow M$, of a
two dimensional oriented, connected,smooth manifold $E$.
 A  tame knot is an embedding $f:S^1 \rightarrow M$ that can be extended
to an embedding of $f: S^1 \times B(0,\epsilon) \rightarrow M$. In
other words tame knots admit a tubular neighborhood embedding. We
will call a manifold \textit{tamely path connected} if for every
two points $P,Q$ of $M$ there is a path $\gamma:[0,1]\rightarrow
M$ with $\gamma(0)=P$, $\gamma(1)=Q$ with the additional property
that for a suitable small disk $B(0,\epsilon)$ the path $\gamma$
can be extended to an embedding $\gamma:B(0,\epsilon) \times [0,1]
\rightarrow M$. It is not clear to the authors whether all path
connected 3 manifolds are tamely path connected. In what follows
we will be concerned only with tamely path connected 3 manifolds.

This is just a small version of the dictionary. More precise
versions can be found in \cite{Mo}, \cite{Si}.

One of the differences between the two theories is that the group
$\mathrm{Cl}(K)$ is always finite while $H_1(M, \mathbb{Z}) =
\mathbb{Z}^r \oplus H_1(M, \mathbb{Z})_{tor}$ is not. Many authors
proposed that the analogue of the class group for arithmetic
topology should be the torsion part  $H_1(M, \mathbb{Z})_{tor}$,
but we think that one advantage, of taking as analogue of the class
group,
 the whole $H_1(M, \mathbb{Z})$ is that $H_1(M,\mathbb{Z})$ is the
Galois group of the \textit{Hilbert manifold} $M^{(1)}$  over $M$,
where  the Hilbert manifold $M^{(1)}$ is the maximal unramified
abelian cover of $M$.

\begin{theorem}[Principal Ideal Theorem for Number Fields.]
Let $K$ be a number field and let $K^{(1)}$ be the Hilbert class
field of $K$. Let $\mathcal{O}_K$, $\mathcal{O}_{K^{(1)}}$ be the
rings of integers of $K$ and $K^{(1)}$ respectively. Let $P$ be a
prime ideal of $\mathcal{O}_{K^{(1)}}$. We consider the prime
ideal
$$\mathcal{O}_{K} \rhd p = P \cap \mathcal{O}_{K}$$ and let
\[
 p\mathcal{O}_{K^{(1)}} = \left (PP_2\ldots P_r \right )^e=\prod_{g\in \mathrm{CL}(K)} g(P)
\]
 be
the decomposition of $p\mathcal{O}_{K^{(1)}}$ in
$\mathcal{O}_{K^{(1)}}$ into prime ideals. The ideal
$p\mathcal{O}_{K^{(1)}}$ is principal in $K^{(1)}$.
\end{theorem}

This theorem was conjectured by Hilbert and the proof was reduced
to a purely group theoretic problem by E. Artin. The group
theoretic question was resolved by Ph. Furtwangler
\cite{Furtwangler}. For a modern account we refer to
\cite[V.12]{Janusz}.

\section{The Principal Ideal Theorem for Knots}
The Hilbert class field in number fields is defined to be the
largest non-ramified abelian extension. Therefore we define the
Hilbert manifold $M^{(1)}$ of $M$ as the universal covering space
$\widetilde{M}$ of $M$ modulo the commutator group $[\pi_1(M),
\pi_1(M)]$:
\[ M^{(1)} = \widetilde{M}/[\pi_1(M), \pi_1(M)]. \]
By definition $M^{(1)}$ is the largest unramified  abelian cover
of the manifold $M$. Moreover, the Galois group of the cover is:
\[G=\mathrm{Gal}(M^{(1)}/M) = \pi_1(M)/[\pi_1(M), \pi_1(M)] = H_1
(M,\mathbb{Z}).\]

Let $L/K$ be a Galois extension of number fields and let
$\mathcal{O}_L,\mathcal{O}_K$ be the corresponding rings of
algebraic integers. In the case of number fields it is known that
every prime ideal $p \lhd \mathcal{O}_K$ gives rise to an ideal $p
\mathcal{O}_L$. This construction is not always possible in the
case of $3$-manifolds. Namely, if $M_1 \rightarrow M$  is a
covering of $3$-manifolds then an arbitrary knot does not
necessarily lift to a knot in $M_1$. Indeed, a knot can be seen as
a path $\gamma:[0,1] \rightarrow M$ so that $\gamma(0)=\gamma(1)$,
and paths do lift to paths $\tilde{\gamma}:[0,1]\rightarrow M_1$,
but in general $\tilde{\gamma}(0)\neq \tilde{\gamma}(1)$.  The
following theorem gives a necessary and sufficient condition for
liftings of maps between topological spaces.

\begin{theorem}\label{4.2}
Let $(Y,y_0)$, $(X,x_0)$ be topological spaces (arcwise connected,
semilocally simply connected), let $p:(X', x_0') \rightarrow
(X,x_0)$ be a topological covering with $p(x_0')=x_0$ and let
$f:(Y,y_0) \rightarrow (X,x_0)$ be a continuous map. Then, there
is a lift $\tilde{f}: Y \rightarrow X'$ of f,

\[
\xymatrix {  & X' \ar[d]^p \\
Y \ar[r]_f \ar@{-->}[ur]^{\tilde{f}} &  X } \]
making the above diagram commutative if and only if
\[f_{\ast} ( \pi_1 (Y,y_0)) \subset p_{\ast}(\pi_1
(X', x_0')),\] where $f_*$, $p_*$ are the induced maps of
fundamental groups.

\end{theorem}
\begin{proof}
\cite[Chapter 5, Proposition 5]{Mas}.
\end{proof}

\begin{pro}\label{knot}
Let $K_1$ be a knot in $M^{(1)}$. Denote by $G(K_1)$ the subgroup of $G$ fixing $K_1$. Consider the link $L=\bigcup_{g
\in G / G(K_1)} gK_1$. Then L is zero in $H_1(M^{(1)},\mathbb{Z})$.
\end{pro}
\begin{proof}
In number theory this theorem is proved by using the transfer map,
but this method can not be applied in our case since $G$ need not
be finite. If $\left | H_1(X,\mathbb{Z}) \right | < \infty  $ then
the classical \cite[V.12]{Janusz} proof applies by just using the
MKR dictionary, i.e. by replacing all the class groups that appear
in the classical proof with the first homology groups. In the
general case we will use the Theorem \ref{4.2}.

Since the diagram
\[
\xymatrix {& K_1  \ar[d]^p  \ar[r]  & M^{(1)} \ar[d]^p\\
 S^1 \ar[r]^f \ar[ur]^{\tilde{f}} & p(K_1)  \ar[r] & M
}
\]
commutes we have that
 \[ f_{\ast} (\pi_1 (S^1)) \subset p_{\ast}(\pi_1(K_1)) \subset
 p_{\ast}(\pi_1 (M^{(1)})) = p_{\ast}([\pi_1(M), \pi_1(M)]),\]
therefore $f_{\ast}(\pi_1(S^1))=0$ as an element in
$H_1(M,\mathbb{Z})$, hence there is a topological (possibly
singular) surface $\phi:E\rightarrow M$  so that
\[f(S^1) = p (K^1) = \partial \phi(E).\]

Moreover the surface $E$ is homotopically trivial therefore
theorem \ref{4.2} implies that  there is a map $\widetilde{\phi}$
making the following diagram commutative:
 \[
 \xymatrix{  \;& M^{(1)}\ar[d]^p \\
 E \ar[r]_{\phi} \ar[ur]^{\widetilde{\phi}} & M },
  \]
with the additional property $\partial \widetilde{\phi}(E) = p^{-1} (\partial \phi(E))=L$.
\end{proof}

Observe that proposition \ref{knot} proves only that there is no topological obstruction for the link
$L$ to be the boundary of a surface. Since we have worked in terms of singular homology the
boundary surface might have singularities or might consist of several components.
We will use the following theorem  known as ``Dehn lemma'' in the literature.
\begin{theorem} \label{DL}
 Let $M$ be a $3$-manifold and $f:D^2 \rightarrow M$ be a map such that for some
neighborhood $A$ of $\partial D^2$ in $D^2$ $f\mid_A$ is an
embedding and $f^{-1} f(A)=A$. Then $f\mid_{\partial D^2}$ extends
to an embedding $g:D^2 \rightarrow M$.
\end{theorem}
\begin{proof}
 \cite[4.1]{Hempel}
\end{proof}

\begin{cor} \label{DLc}
 If a tame knot is the boundary of  a topological and possibly singular surface then the knot
is the boundary of an embedded surface.
\end{cor}
\begin{proof}
 Using the embedding of a  tubular  neighborhood of the knot we can construct a nonsingular
collar around  the boundary of the topological surface and the desired result follows by theorem \ref{DL}.
\end{proof}
\begin{pro}
 Let $L$ be a link in $M$ that is a homologically trivial. Then there is a family of tame knots  $K_\epsilon$ in M with  $\epsilon>0$,
that are boundaries of  embedded surfaces $E_\epsilon$ so that
$\lim_{\epsilon \rightarrow 0} K_\epsilon=L$ and $E=\lim_{\epsilon
\rightarrow 0} E_\epsilon$ is an embedded surface with $\partial
E=L$.
\end{pro}
\begin{proof}
 We will consider the case of a link with two components. Let $L=K_1 \cup K_2$, where $K_i$ is given by the embedding  $f_i:S^1 \rightarrow M$, a tame knot.
 The passage from two components to $n>2$ components follows by induction.
Select two points $P_1,Q_1$ on $f_1(S^1)$ and two points $P_2,Q_2$
on $f_2(S^1)$ so that $d(P_i,Q_i)=\epsilon$. The embedding  $f_i$
can be given as the union of two curves, namely
$\gamma_i:[0,1]\rightarrow M$, $\delta_i:[0,1] \rightarrow M$, so
that $\gamma_i(0)=\delta_i(1)=P_i$, $\gamma_i(1)=\delta_i(0)=Q_i$.
This means that the ``small'' curve is the curve $\delta_i$.

Since the manifold $M$ is tamely path connected we can find two
paths $\alpha,\beta:[0,1]\rightarrow M$ so that
$\alpha(0)=P_1,\alpha(1)=Q_2$, $\beta(0)=P_2,\beta(1)=Q_1$, that
are close enough so that the rectangle $\alpha  (-\delta_2) \beta
 (-\delta_1)$ is homotopically trivial.
Let $I=[0,1]\subset \mathbb{R}$.
Every path in $M$, {\em i.e.} every function $f:I\rightarrow M$, defines a cycle in $H_1(M,\mathbb{Z})$.
We will abuse the notation and we will denote by $f(I)$ the homology class of the path $f(I)$.
We compute in  $H_1(M,\mathbb{Z})$:
\[
 0=f_1(S_1) + f_2(S_1)=\gamma_1(I)+\gamma_2(I)+\delta_1(I)+\delta_2(I)+0=\]
\[=
\gamma_1(I)+\gamma_2(I)+\delta_1(I)+\delta_2(I)+ \alpha(I) -\delta_2(I) + \beta(I) -\delta_1(I)=
\]
\[
 =\gamma_1(I)+\alpha(I)+\gamma_2(I)+\beta(I).
\]
This means that the tame knot $\gamma_1 \alpha  \gamma_2 \beta$ is
the boundary of
 a topological surface, and by Corollary \ref{DLc} it is the boundary of an embedded surface $E_\epsilon$.

Choose an orientation on $E_\epsilon$ so that on $P\in \partial
E_\epsilon$ one vector of the oriented basis of $T_P E_\epsilon$
is the tangent vector of the curves $\partial E_\epsilon$ and the
other one is pointing inwards of $E$. We will denote the second
vector by $N_P$. Moreover, we choose the same orientation on all
surfaces $E_\epsilon$ in the same way, {\em i.e.} the induced
orientation on the common curves of the boundary is the same.

We would like to take the limit surface for $\epsilon \rightarrow
0$. For this we have to distinguish the following two cases: In
the first case the direction of decreasing the distance $\epsilon$
is the opposite of $N_P$ and the limiting procedure makes the
rectangle $\alpha \cdot (-\delta_2) \cdot \beta  \cdot
(-\delta_1)$ thinner and eventually it eliminates it. In this case
the elimination of the above mentioned rectangle glues two parts
of the surface $E_\epsilon$ together. The limit  $\epsilon
\rightarrow 0$ gives us an embedded  surface $E$ that is the
boundary of our initial link $L$. Indeed by taking the limit the
paths $\alpha(I), \beta(I)$ are identified, and this
identification can be done in a smooth manner.

In the second case the direction of decreasing the distance
$\epsilon$ is the same with  $N_P$. This means that by taking the
limit $\epsilon \rightarrow 0$ we don't glue two parts of the
boundary of the surface $E_\epsilon$ but we make the rectangle
$\alpha  (-\delta_2)  \beta   (-\delta_1)$ thinner and after
eliminating it we cut the surface in two pieces. Still the limit
$\epsilon \rightarrow 0$ gives us two embedded surfaces $E,E'$
that are the boundaries of our initial link components $K_1,K_2$.
We can arrive at one embedded surface in the following way: We
cut two disks $D_1,D_2$  of the interiors of $E$ and $E'$
  and glue together them together along a tubular path $T$ so that $\partial T=D_1\cup D_2$.
\end{proof}

As a corollary of the principal ideal theorem for knots we state
the following:
\begin{theorem}[\textit{Seifert}]
 Every link in a simply connected $3$ manifold  is the boundary of an embedded surface.
\end{theorem}
\begin{proof}
Let $M$ be simply connected. The Hilbert manifold of $M$ coincides with $M$ and the result follows.
\end{proof}

\end{document}